\documentstyle[12pt,amssymb,amstex,times]{amsart}

\numberwithin{equation}{section}

\def\P{\mathcal{P}}
\def\a{\alpha}

\def\d{\delta}  

\def\e{\varepsilon}
\def\l{\lambda}

\def\r{\rho}
\def\s{\sigma} 
\def\t{\tau}

\def\id{{\bf 1}\!\!{\rm I}}
\newtheorem{thm}{Theorem}[section]
\newtheorem{lem}[thm]{Lemma}
\newtheorem{cor}[thm]{Corollary}
\newtheorem{prop}[thm]{Proposition}
\newtheorem{defin}[thm]{Definition}
\newtheorem{rem}[thm]{Remark}

\begin{document}

\title[]
{On Limit theorems in $JW$- algebras}

\author{Abdusalom Karimov}
\address{Abdusalom Karimov\\
Department of Mathematics\\
Tashkent Institute of Textile and Light Industry\\
Tashkent, 100100, Uzbekistan} \email{\tt karimov57@@rambler.ru}

\author{Farrukh Mukhamedov}
\address{Farrukh Mukhamedov\\
 Department of Computational \& Theoretical Sciences\\
Faculty of Science, International Islamic University Malaysia\\
P.O. Box, 141, 25710, Kuantan\\
Pahang, Malaysia} \email{{\tt far75m@@yandex.ru}, {\tt
farrukh\_m@@iiu.edu.my}}

\begin{abstract}
In the present paper, we study bundle convergence in $JW$- algebra
and  prove certain ergodic theorems with respect to such
convergence. Moreover, conditional expectations of reversible
$JW$-algebras are considered. Using such expectations, the
convergence of supermartingales in such is established.

\end{abstract}
\maketitle

\section{Introduction.}
The concept of almost everywhere convergence and its different
variants for sequences in von Neumann algebra were studied by many
authors (see for example, \cite{J, L, P}). There were proved many
limit and ergodic theorems with respect to almost everywhere
convergence in such algebras with faithful normal state \cite{Ba, G,
Y}. In \cite{HJP}, a notion of bundle convergence in a von Neumann
algebra was firstly defined which coincides with usual almost
everywhere convergence in the case of commutative algebra
$L_{\infty}$. Certain limit theorems with respect to such a
convergence were obtained there. Other several results concerning
the bundle convergence in a von Neumann algebra and in its
$L_2$-space were studied in \cite{BF1,BF2,S}.

On the other hand, in most mathematical formulations of the
foundations of quantum mechanics, the bounded observables of a
physical system are identified with a real linear space, $L$, of
bounded self-adjoint operators on a Hilbert space $H$. Those
bounded observables which correspond to the projections in $L$
form a complete orthomodular lattice, $P$, otherwise known as the
lattice of the quantum logic of the physical system. For the
self-adjoint operators $x$ and $y$ on $H$ their Jordan product is
defined by $x\circ y=(xy+yx)/2=(x +y)^2-x^2-y^2$. Thus $x\circ y$
is self-adjoint so it is reasonable to assume that $L$ is a Jordan
algebra of self-adjoint operators on $H$ which is closed in the
weak operator topology. Hence L is a $JW$-algebra. It is known
that the $JW$- algebra is a real non-associative analog of a von
Neumann algebra that firstly studied by Topping \cite{T}. He had
extended many results from the theory of von Neumann algebras on
$JW$- algebras. Particulary, in \cite{A2}, a problem of extension
of states and traces from $JW$- algebra to its enveloping von
Neumann algebra was solved. In \cite{A3}, concepts of convergence
in measure, almost uniformly convergence and convergence
$s$-almost uniformly in $JW$-algebras were introduced, and
relations between them were investigated as well. Such kinds of
convergence are used to prove various ergodic theorems for Markov
operators \cite{A1,A4,A51,K,KM1}. Asymptotic behavior of positive
contractions of Jordan algebras has been studied in
\cite{MTA,MTA2}. We refer the reader to the books \cite{A5,ARU,HS}
for the theory of $JW$-algebras.

The main purposes of this paper is to extend a concept of the bundle
convergence from von Neumann algebras to $JW$-algebras and prove
certain limit theorems. The paper is organized as follows. In
section 2, we recall some well-known facts and basic definitions
from the theory of Jordan and von Neumann algebras. In section 3,
the bundle convergence and its properties are studied. It is proved
that the bundle limit of a sequence of uncorrelated operators is
proportional to the identity operator. In the last section 4,
certain properties of conditional expectations of a reversible $JW$-
algebra $A$ with a faithful normal trace are studied, and a
martingale convergence theorem is proved as well. Here, we apply a
method passing from Jordan algebra to the corresponding enveloping
von Neumann algebra \cite{K}. Note that supermartingales and
martingales in a von Neumann algebra setting were investigated in
\cite{Br,G,J}.

\section{Preliminaries}

Throughout the paper $H$ denotes a complex Hilbert space, $B(H)$
denotes  the algebra of all bounded linear operators on $H$.

Recall that a {\it $JW$-algebra} is a real linear space of
self-adjoint operators from $B(H)$ which is closed under the Jordan
product $a\circ b={1\over 2}(ab+ba)$ and also closed in the weak
operator topology. Here the sign $ab$ denotes the usual operator
multiplication of operators $a$ and $b$ taken from $B(H).$ A $JW$-
algebra $A$ is said to be {\it reversible} if $a_1a_2\cdots
a_n+a_na_{n-1}\cdots a_1\in A$, whenever $a_1, a_2, \dots , a_n\in
A$. Examples of non-reversible $JW$-algebras are spin factors which
are described in \cite{T,ARU}.

Recall that a real $*$-algebra $R$ in $B(H)$ is called a {\it real
$W^*$-algebra} if it is closed in the weak operator topology and
satisfies the conditions $R\cap i R=\{0\}$, $\id\in R$. It is
obvious that if $R$ is a real or complex $W^*$-algebra, then its
selfadjoint part $R_{sa} =\{x\in R : x^* = x\}$ forms a reversible
$JW$-algebra.

Given an arbitrary $JW$-algebra $A$ let $R(A)$ denote the weakly
closed real $*$-algebra in $B(H)$ generated by $A$, and let $W(A)$
denote the $W^*$- algebra (complex) generated by $A$.

\begin{thm}\label{2.4} \cite{HS} Let $A$ be a reversible  $JW$-algebra, then
the following assertions hold
\begin{enumerate}
\item[(i)] $R(A)_{sa}=A$ and $W(A)=R(A)+iR(A)$;

\item[(ii)] $\| a+ib\| \geq \max \{\| a \|, \| b \| \} $ for every  $a,b\in
R(A) $;

\item[(iii)] if $a+ib\geq0$ for $a,b\in R(A) $, then $a \geq 0$.
\end{enumerate}
\end{thm}

Throughout the paper we always assume that a $JW$-algebra $A$ is
reversible, therefore we do not stress on it, if it is not
necessary.

Let $A$ be a reversible $JW$- algebra and $\r$ be a faithful
normal (f.n.)  state (resp.  faithful normal semifinite (f.n.s.)
trace) on $A$. Then $\r$ can be extended to a f.n. state $\r_1$
(resp. f.n.s. trace) on the $W^*$-algebra $W(A)$. Namely, for
every $x\in W(A)^+$ ($W^+$ means the positive part of $W$), one
has $x=a+ib$, where $a,b\in R(A)$, $a\in A^+$, $b^*=-b$
(skew-symmetric element). Then, one puts $\r_1(x)=\r(a)$ (see
\cite{ARU}, Thm. 1.2.9 for more details).

Let $A$ be  a  $JW$- algebra with a f.n.s. trace $\t$ and $\t_1$ be
its  extension to $W(A)$. Set ${\frak N}_{\t}=\{x\in A:  \t
(|x|)<\infty\}$, ${\frak N}_{\t_1}=\{x\in W(A): \t_1(|x|)<\infty\}$.
Completion of ${\frak N}_{\t}$ (resp. ${\frak N}_{\t_1}$) w.r.t. the
norm $\|x\|_1=\t (|x|)$ \ $x\in{\frak N}_{\t}$  (resp. $\|x\|_1=\t_1
(|x|)$ \ $x\in{\frak N}_{\t_1}$) is denoted by $L_1(A, \t)$ (resp.
$L_1(W(A), \t_1)$). It is obvious that $L_1(A, \t)\subset L_1(W(A),
\t_1).$

Let $W$  be a von Neumann algebra with a f.n. state $\r$. By $Proj
W$ we denote the set of all projections in $W$.

We recall some facts about almost uniformly convergence.

A sequence $\{x_n\}\subset W$ is said to be {\it almost uniformly
convergent} to $x\in W$ ($x_n \buildrel {a. u.} \over
\longrightarrow x$) if, for each $\e >0$, there exists $p\in Proj W$
with $\r(p^{\bot})<\e$ such that $\|(x_n-x)p\|\to 0$ as
$n\to\infty.$

Further, we need the following

\begin{lem}\label{2.6} \cite{J} For a uniformly bounded sequence, almost uniformly convergence
implies strong convergence.

\end{lem}

Now we give some necessary definition and results from \cite{HJP}.

Suppose that $\{D_m\}\subset W^+$ with $\sum\limits_{m = 1}^\infty
{\r(D_m) }<\infty$. {\it The bundle} (determined by the sequence
$\{D_m\}$) is the set
$$
P_{(D_m )} =\bigg\{p \in Proj W: p\neq0, \ \ \sup\limits_m \bigg\|
{p\bigg(\sum\limits_{k = 1}^m {D_k }\bigg)p} \bigg\| < \infty \
\textrm{and} \ \|pD_mp\| \to 0 \ \textrm{as} \ m\to \infty \bigg\}.
$$

Let $x_n, x\in W$ $(n=1,2,\dots)$. We say that $\{x_n\}$ is {\it
bundle convergent} to $x$ ($x_n \buildrel {b, W} \over
\longrightarrow x$) if there exists a bundle $P_{(D_m)}$  with
$D_m\in W^+, \sum\limits_{m = 1}^\infty {\r(D_m) }<\infty$ such that
$p\in P_{(D_m)}$ implies $\|(x_n-x)p\|\to 0.$

\begin{lem}\label{2.8} \cite{HJP} Let $0<\e<1/16,$ $D_m\in W^+$ \
$(m=1,2,\dots)$ and $$\sum\limits_{m = 1}^\infty {\r(D_m) }<\e.$$
Then there exists $p\in Proj\ W$  such that
$$\r(p^\bot)<\e^{1/4},  \   \   \
 \left\|{p\left({\sum\limits_{k = 1}^m {D_k }}\right)p}\right\| < 4\e^{1/2}, m=1,2,\dots$$

\end{lem}

\section{The bundle convergence in $JW$- algebras}

In this section, we study bundle convergence and its properties in $JW$-algebras.
We shall prove that the bundle limit of a sequence of uncorrelated
operators is proportional to the identity operator. Throughout the
paper $A$ denotes a reversible $JW$- algebra.

 Let $A$ be a
$JW$- algebra with a f.n. state $\r$ and $\r_1$ be its extension to
the enveloping von Neumann algebra $W(A)$.

Suppose that $\{ {D_m }\} \subset A^ + $ with $\sum\limits_{m =
1}^\infty {\r (D_m )}  < \infty .$  {\it The bundle} is the set
$${\P}_{(D_m )}  =\bigg\{p \in Proj A: p\neq0,  \ \sup\limits_m
\bigg\| {p(\sum\limits_{k = 1}^m {D_k } )p} \bigg\| < \infty \
\textrm{and} \  \|p(D_m)p\| \to 0 \ \textrm{as} \ m\to \infty
\bigg\}.$$

\begin{rem}\label{3.2} Let $\{D_m\}^{\infty}_{m=1}\subset A^+$,
$\sum\limits_{m = 1}^\infty \r (D_m)< \infty$. Then
${\P}_{(\{D_m\}^\infty _{m=1})}={\P}_{(\{D_m\}^\infty _{m=k})}$ and
in the definition of bundle we may suppose that \ $\sum\limits_{m =
1}^\infty \r (D_m)< \e$, for some positive number $\e$.

\end{rem}

 Let $x_n$, $x\in A$ $(n=1,2,\dots).$
We say that $\{x_n\}$ is {\it bundle convergent to} $x$ ($x_n
\buildrel {b, A} \over \longrightarrow x$) if there exists a bundle
${\P}_{(D_m)}$ such that $p\in {\P}_{(D_m)}$ implies
$\|p(x_n-x)^2p\|\to 0.$

Clearly, an intersection of two bundles is a bundle as well.
Consequently, the bundle convergence in $A$ is additive, and the
bundle limit in $A$ is unique. In the case when $A$ is a
self-adjoint part of some von Neumann algebra,  bundle convergence
in Jordan algebra coincides with convergence in a von Neumann
algebra setting.

Similarly, one can define that a sequence $\{x_n\}\subset A$ is said
to be {\it almost uniformly convergent} to $x\in A$ ($x_n \buildrel
{a. u.} \over \longrightarrow x$) if, for each $\e >0$, there exists
$p\in Proj A$ with $\r(p^{\bot})<\e$ such that $\|p(x_n-x)^2p\|\to
0$ as $n\to\infty$ \cite{A3}.

In what follows we need the following auxiliary result.

\begin{lem}\label{pdp} Let $\s,\d>0$, $D\in A^+$ and suppose
$p\in Proj W(A)$ with $\r_1(p^{\bot})<\s$ such that $\|pDp\|<\d$.
Then there exists a projection $E\in Proj A$ with $\r(E^{\bot})<2\s$
such that $\|EDE\|<4\d$
\end{lem}
\begin{pf} Suppose $p=x+iy$, \ $x,y\in R(A)$,  one finds $0\leq x\leq \id$ and let $x =
\int\limits_0^1 {\l de_\l} $ be the spectral resolution of $x$. Put
$$b =\int\limits_{{1/ 2}}^1 {\l ^{ - 1} de_\l }.$$
Then
$$
E=xb = \int\limits_0^1 {\l de_\l} \int\limits_ {{1 / 2}}^1 {\l ^{
- 1} de_\l  } =s(x) - e_{1 / 2}.$$ ($s(x)$- means the support
projector of $x$) and
\begin{eqnarray*} 2(s(x)-x)&=& 2\int\limits_0^1 {(1-\l) de_\l}\geq
2\int\limits_0^{1 /2} {(1-\l) de_\l}\\&>&\int\limits_0^{1/2} {
de_\l}= e_{1/2}.\end{eqnarray*} Since  \  $s(x)^{\bot}\leq
2s(x)^{\bot}$ \ then  \ $\id -E=s(x)^{\bot}+e_{1/2}\leq
2(s(x)^{\bot}+(s(x)-x))=2(\id -x)$ \  and  \  $\r (\id- E)\leq2\r
(\id - x)= 2\r_1(\id-x-iy)=2\r_1(\id-p)<2\s$. This means that
$E\neq0$.

Then the inequality (ii) of Theorem \ref{2.4} with commutativity of
$b$ and $x$ implies that
\begin{eqnarray*}
\|EDE\|&=&\|xbDxb\|\leq
\|b\|^2\|xDx\|\\
&=&4\|\sqrt{D}x\|^2 \leq 4\|\sqrt{D}(x+iy)\|^2\\
&=4\|pDp\|,
\end{eqnarray*} which means $\|EDE\| <4\d$. This completes the proof.
\end{pf}

Using Lemma \ref{2.8} and \ref{pdp} one gets the following

\begin{lem}\label{3.4} Let $0<\e<1/16,$ $D_m\in A^+$ \
$(m=1,2,\dots)$ and $$\sum\limits_{m = 1}^\infty {\r(D_m) }<\e.$$
Then there exists $E\in Proj A$  such that
$$\r(E^{\bot})<2\e^{1/4},  \   \   \
\left\|{E\left({\sum\limits_{k = 1}^m {D_k }}\right)E}\right\| <
4\e^{1/2}, m=1,2,\dots$$
\end{lem}

\begin{cor}\label{3.5}
For each bundle ${\P}_{(D_m)}$ (with $D_m\in A^+$ and
$\sum\limits^\infty_{m=1}\r(D_m)<\infty$) and for each $\e>0$ there
exists $p\in {\P}_{(D_m)}$ such that $\r(p^\bot)<2\e.$
\end{cor}

\begin{pf} According to Remark \ref{3.2} suppose that $\sum\limits_{m = 1}^\infty
\r (D_m)<\e $. Let $0<\a_m\nearrow \infty$ be a sequence of
numbers such that $\sum\limits^\infty_{m=1}\a_m\r(D_m)<\e$. Put
$B_m=\a_mD_m$ $m=1,2,\dots$ and applying Lemma \ref{3.4} to $B_m$
and $\e^4$ we get the existence of $p\in {\P}_{(B_m)}$ such that
$\r(p^\bot)< 2 \e$, $\left\|{p\left({\sum\limits_{k = 1}^m \a_k
{D_k }}\right)p}\right\| < 4\e^2<\infty$, but
$\|pD_mp\|\leq\a_m^{-1}\left\|\sum\limits_{k = 1}^m p(\a_k {D_k
)p}\right\| \leq 4\e^2 \a^{-1}_m$ and $\|pD_mp\|\to0$, so  $p\in
{\P}_{(D_m)}.$
\end{pf}

Obviously, this corollary implies the following

\begin{prop}\label{3.6}  If $x_n \buildrel {b, A} \over \longrightarrow x$
then  $x_n \buildrel {a. u.} \over \longrightarrow x$.
\end{prop}

\begin{rem}\label{3.7} Let $\{x_n\}$ be  uniformly bounded   and $x_n \buildrel{a.u.}
\over \longrightarrow x$  in  $A$, $x\in A,$ then $x_n \to x$  in
the strong    operator topology.
\end{rem}

Indeed, let $x_n,$  $x\in A$, then we obviously have   $x_n
\buildrel{a.u.} \over \longrightarrow x$ in $W(A).$  From Lemma
\ref{2.6} we derive that $x_n \to x$ in the strong operator
topology, i.e. for every $\xi\in H$ $\|x_n\xi-x\xi\|_H\to 0$.

\begin{prop}\label{3.8} Let $\{x_n\}\subset A$ and
$\sum\limits_{n=1}^\infty\r(|x_n|^2)<\infty$. Then $x_n\buildrel
{b,A} \over\longrightarrow 0$.

\begin{pf} The sequence $D_m=|x_m|^2$, \ $m=1, 2, \dots$ defines a
bundle ${\P}_{(D_m)}$. Let $p\in {\P}_{(D_m)}$. Then $\|pD_mp\|\to
0$ and $\|px_n^2p\|=\|x_np\|^2=\|p|x_n|^2p\|=\|pD_mp\|$.
\end{pf}
\end{prop}

\begin{thm}\label{bc-au} Let $A$ be a reversible $JW$-algebra with  f.n. state
$\rho$, and $\rho_1$ be its extension to the
enveloping von Neumann algebra $W(A)$. If a sequence $(x_n)\subset
A$ is bundle convergent in $(W(A);\rho_1)$, then it is bundle
convergent in $(A;\rho)$.
\end{thm}

\begin{pf} Suppose that $x_n \buildrel b,W(A) \over
\longrightarrow x$. Then there is a sequence $\{D_m\}\subset W(A)^+$
\ $m=1,2,\dots$ with $\sum\limits_{m = 1}^\infty {\rho_1 (D_m )} <
\infty $ and the corresponding bundle
$$
P_{(D_m )}  =\bigg\{p \in Proj W(A),  p\neq0:  \sup_m \bigg\|
{p(\sum\limits_{k = 1}^m {D_k } )p} \bigg\| < \infty \
\textrm{and} \ \|pD_mp\| \to 0 \ \textrm{as} \ m\to \infty \bigg\}
$$
such that  for each $p\in P_{(D_m)}$ one has $\|(x_n-x)p\|\to 0.$
Due to Theorem \ref{2.4} one has $D_m =K_m+iL_m\geq 0$, and
$K_m\in A^+$. Hence, according to the definition of $\rho_1$ we
derive $\rho(K_m)=\rho_1(K_m+iL_m).$ Therefore, $\sum\limits_{m =
1}^\infty {\rho (K_m )} < \infty .$ By $P_{(K_m)}$ we denote the
bundle in $W(A)$, generated by the sequence $\{K_m\}$.  Let $p\in
P_{(D_m)}$, then (ii) of Theorem \ref{2.4} yields that
$\|pK_mp\|\leq \|pK_mp+ipL_mp\|=\|pD_mp\|\to 0$. Similarly, one
can establish that $\mathop {\sup }\limits_m \left\|
{p(\sum\limits_{i = 1}^m {K_i } )p} \right\| < \infty.$ Hence,
$P_{(D_m)}\subset P_{(K_m)}$, and so $P_{(K_m)}$ is nonempty.

Now we show that $P_{(K_m)}\subset P_{(D_m)} $. Let us consider a
mapping $\eta:W(A)\to W(A)$ defined by $\eta(x+ iy)=x-iy$, which is
a $*-$ anti-automorphism of $W(A)$ \cite{HS}. If $x+iy=a^2\geq 0$
for some $a\in W_{sa}(A)$, then
$x-iy=\eta(x+iy)=\eta(a^2)=(\eta(a))^2\geq 0$, here $x^*=x$,
$y=-y^*$. This yields $x-iy\geq 0$, and hence we have $x+iy\leq 2x$.
Now applying the last inequality to $K_m+iL_m$ we infer that
$D_m\leq 2K_m$. Assume that $p\in P_{(K_m)}$, then  $\|pK_mp\|\to 0$
 as well as $2\|pK_mp\|\to 0.$ Consequently, we find  $\|pD_mp\|\leq 2\|pK_mp\|\to 0$
which implies that $p\in P_{(D_m)}$. So, $P_{(D_m)}\ = P_{(K_m)}$.

Put ${\P}= P_{(K_m)}\cap A$, which is the bundle  in $A$ generated
by $\{K_m\}.$ Then ${\P}\subset P_{D_m}$ and the proof is
finished.

  \end{pf}

As a consequence of the last result we have the following result.

\begin{thm}\label{3.10}  Let $A$ be a reversible $JW$-algebra with  f.n. state
$\rho$, and $\rho_1$ be its extension to the
enveloping von Neumann algebra $W(A)$. Suppose that $\{x_n\}$ is a sequence of pairwise
orthogonal operators in $A$ (i.e. $\r(x_i\circ x_j)=0$ for $i\neq
j$). If
$$\sum\limits_{n=1}^\infty \bigg(\frac{\log_2(n+1)}{n}\bigg)^2\r(|x_n|^2)<\infty
$$  then $\frac{1}{n}\sum\limits_{j=1}^nx_j\buildrel {b,A}
\over\longrightarrow 0$.
\end{thm}

\begin{pf} Let us first note that if $a, b\in A$ are orthogonal in $A$ with
respect to   $\r$, i.e. $\r(a\circ b)=0$, then for the extended
state $\r_1$ we have $\r_1(ab + ba)= 2\r(a\circ b)=0$. Taking into
account skew-symmetricity of $ab-ba$ and the definition of the
extended state $\r_1$, i.e. the fact $\r_1(x)=0$ whenever $x$ is
skew-symmetric, we got $\r_1(ab - ba)=0$, which with $\r_1(ab +
ba)=0$ implies that $\r_1(ab)=0$, i.e. $a$ and $b$ are orthogonal in
$W(A)$. Hence, the elements $x_n$, $n\in\mathbb{N}$ are pairwise
orthogonal in $W(A)$, so due to Theorem 4.2 \cite{HJP} one finds
$\frac{1}{n}\sum\limits_{j=1}^nx_j\buildrel {b,W(A)}
\over\longrightarrow 0$. Consequently, by means of Theorem
\ref{bc-au} we got the desired statement.
\end{pf}

We have the following

\begin{lem}\label{3.11} Let $x_n,$ $x\in A$ $(n=1,2,\dots)$, \ and \
$x_n \buildrel b,A \over \longrightarrow x,$ \ \ then $\frac{1}{n}
\sum\limits_{k = 1}^n {x_k \buildrel b,A \over
 \longrightarrow x}. $

\end{lem}

Recall that {\it the covariance} of two operators $a,b\in A$ is
defined by $cov(a,b)=\r(a\circ b)-\r(a)\r(b),$ {\it the variance} of
$a$ is defined by $var(a)=cov(a,a)=\r(a\circ a)-|\r(a)|^2=
\r(a^2)-|\r(a)|^2.$ If $cov(a,b)=0$,  then $a$ and $b$ are called
{\it uncorrelated}.

\begin{thm}\label{3.13} Let $A$ be a reversible $JW$-algebra, acting on Hilbert space
$H$, with  f.n. state $\rho$. Suppose that  $\{x_n\}$ is  a sequence of uncorrelated operators in $A$,
and the following conditions are satisfied:
\begin{enumerate}
\item[(i)] $x_n \buildrel {b, A} \over \longrightarrow x$ \ \ $(x\in
A),$

\item[(ii)] $\sum\limits_{n = 1}^\infty
n^{-2}var(x_n)\log^2(n+1)<\infty.$
\end{enumerate}
Then there exists  a  complex number $c$ such that $x=c \id$, where, as before,
$\id$ is the identity in $A.$
\end{thm}

\begin{pf} Put $y_n=x_n-\r(x_n)\id,$ then $\r(y_m\circ y_n)=0$,  $m\neq n$, and
$\r(y_n^2)=var(x_n),$  \ $m,n=1,2,\dots$ Consequently, Theorem
\ref{3.10} implies
\begin{equation}\label{(3)}
\frac{1}{n}\sum\limits_{k = 1}^n {y_k \buildrel
b,A \over \longrightarrow 0}.
\end{equation}
Let us  put $M_n= \frac{1}{n}\sum\limits_{k = 1}^n x_k  $  and $\l_n=\r(M_n)$,
 for every  $n\geq1$. It then follows  \eqref{(3)}
that $M_n - \l_n\id \buildrel b,A \over \longrightarrow 0$. From (i)
with Lemma \ref{3.11} one finds $M_n\buildrel b,A \over
\longrightarrow x$. Now the additivity of bundle convergence yields
$\l_n\id \buildrel {b, A} \over \longrightarrow x$, i.e. there
exists a bundle ${\P}$ such that $p\in{\P}$, \
$\|p(\l_n\id-x)^2p\|\to 0.$ This means that the sequence $\{\l_np\}$
is uniformly bounded in the norm $\|\cdot \|$, which implies that
the sequence of complex numbers $\{\l_n\}$ is bounded as well. Then
from Remark \ref{3.7} we have $\l_n\id\to x$ in the strong operator
topology, i.e. for every $\xi\in H$ one has  $\|\l_n\xi-x\xi\|_H\to
0.$

Since $\{\l_n\}$ is bounded sequence, then there exists a subsequence $\{\l_{n_k}\}$
such that $\l_{n_k}\to \l$ for some $\l\in\mathbb{C}$. This yields that
$\l_{n_k}\id\to\l\id$ in the strong operator topology. By the strong
convergence of $\l_{n_k}\id$ to $x$, so we derive  $x=\l\id$. This completes the proof.
\end{pf}

Note that if instead of a state $\r$ one takes f.n. finte trace,
then (i) assumption could be replaced with $x_n \buildrel {a.u.}
\over \longrightarrow x$. Note that, in general, such a result might
not valid, since a.u. convergence is not additive (see \cite{Pa}).

\section{Conditional expectation in $JW$- algebras.}

In this section, certain properties of conditional expectations of
reversible $JW$- algebras are studied. Using such properties prove
a main result of this section, formulated in Theorem \ref{4.8},
which is a Jordan analog of the following result:

\begin{thm}\label{4.1} \cite{Br} Let $W$ be a von Neumann algebra with a faithful normal
semifinite trace $\t$. Let $\{x_{\a}\}$ be a supermartingale in
$L_1(W, \t).$ If $\{x_\a\}$ is weakly relatively compact, then
there is an
 $x \in L_1(W, \t)$ such that $x_\a\to x$ in
$L_1(W, \t).$
\end{thm}

Let $A$ be a reversible $JW$- algebra with a f.n.s. trace $\t$. Let
$A_1$ be its $JW$- subalgebra containing the identity operator
$\id$.

Recall that a linear mapping $\phi:A\to A_1$ is called {\it a
conditional expectation} with respect to a $JW$- subalgebra $A_1$ if
the following conditions are satisfied:
\begin{enumerate}
\item[(i)] $\phi(\id)=\id;$
\item[(ii)] If $x\geq 0$, then $\phi(x)\geq 0$;
\item[(iii)] $\phi(xy)=\phi(x)y$ for $x\in A,$ $y\in A_1.$
\end{enumerate}
 Let $\widetilde\t:=\t\upharpoonright{A_1}$ be the restriction of the trace $\t$
to $A_1$ such that $\widetilde\t$ is also semifinite. Then the space
$L_1(A_1,\widetilde\t)$ of integrable operators w.r.t. $(A_1,
\widetilde\t)$, is a subspace of $L_1(A,\t)$.

\begin{thm}\label{4.3} \cite{Ber} Let $A$ be a reversible $JW$- algebra
with a f.n.s trace $\t$ and  $A_1$ be its $JW$- subalgebra with
$\id$. Let $\widetilde\t=\t\upharpoonright{A_1}$ be the
restriction of the trace $\t$ to $A_1$ such that $\widetilde\t$ is
also semifinite. Then  there exists a unique positive linear
mapping $E(\cdot/A_1):A\to A_1$ satisfying the condition $$\t
(E(a/A_1)b)=\t(ab)$$ for $a\in A,$ $b\in L_1(A_1,\widetilde\t)$,
and the mapping $E(\cdot/A_1)$ is a conditional expectation with
respect to $A_1$.
\end{thm}
Such a conditional expectation is called $\t-$ invariant. Note
that when trace is finite an analogous result has been proved in
\cite{A4}. In what follows $E(\cdot/A_1)$ and its extension to
$L_1(A, \t)$ is called a {\it conditional expectation}  w.r.t.
$A_1$ \cite{Ber}.

We note that if $W(A)$ is the enveloping von Neumann algebra of a
$JW$-algebra $A$, then the existence of conditional expectations
from $W(A)$ onto $A$ has been proved in \cite{HaS}.

Let $A$ be a reversible $JW$- algebra with f.n.s. trace $\t$ and
$\t_1$ be its extension to $W(A)$. Further, we suppose that the
restriction of $\t$ to any considered subalgebras is semifinite.
If $B$ is a reversible $JW$- subalgebra of $A$, then we denote a
conditional expectation with respect to $B$ by $E(\cdot/B)$.
Similarly, a conditional expectation from $W(A)$ to a subalgebra
$W(B)$ is denoted by $\widetilde E(\cdot/W(B))$.

So let $\widetilde E(\cdot/W(B))$ be a $\t-$ invariant conditional
expectation, i.e. \begin{equation}\label{(2)} \t_1(\widetilde
E(a/W(B))x) =\t_1(ax)
\end{equation}
   which is uniqe (see \cite{U} )

\begin{lem}\label{4.4} Let $z\in {\frak N}_{\t_1}\cap R(A)$, and
$z^*=-z$, then $\t_1(z)=0$.
\end{lem}

\begin{pf} Let $z\in {\frak N}_{\t_1}\cap R(A)$, and
$z^*=-z$, so $(iz)^*=iz$, $iz\in {\frak N}_{\t_1}.$ From
definition of ${\frak N}_{\t_1}$ we have  $|z|\in {\frak
N}_{\t_1}$. Hence $|z|\in {\frak N}_{\t}.$
 Since $|z|+iz \geq 0,$ then
$\t_1(|z|+iz)=\t_1(|z|).$ According to the linearity of $\t_1$ on
${\frak N}_{\t_1},$  \  $\t_1(|z|+iz)=\t_1(|z|)+i\t_1(z),$ so
$\t_1(|z|)=\t_1(|z|)+i\t_1(z)$ and $\t_1(z)=0.$
\end{pf}

It is naturally to ask: how the restriction of $\widetilde
E(\cdot/W(B))$ to $A$ is related to $E(\cdot/B)$? Next result
answers to such a question.

\begin{thm}\label{4.5} The restriction of $\t-$ invariant conditional
expectation $\widetilde E(\cdot/W(B))$ to $A$ is equal to
$E(\cdot/B),$ here as before $E(\cdot/B)$ is $\t-$ invariant
conditional expectation.

\end{thm}

\begin{pf} It is sufficient to prove that for every $a\in
A $ \begin{equation}\label{(2)}
 \t_1(ax)=\t_1(E(a/B)x)
\end{equation}
 holds for any $x\in
L_1(W(B),\t_1)$, where
$L_1(W(B),\t_1)=L_1(R(B),\t_1)+iL_1(R(B),\t_1),$ and here
$L_1(R(B),\t_1)$ is the completion of  $R(B)\cap {\frak N}_{\t_1}$
with respect to $L_1$-norm. We firstly prove \eqref{(2)} for $x\in
L_1(R(B),\t_1).$ A functional $\t_1(h\cdot)$ is $L_1$-continuous on
$L_1(W(B),\t_1)$, for any $h\in W(B),$ therefore it is enough to
prove \eqref{(2)} for any $x$ taken from $R(B)\cap{\frak N}_{\t_1}.$

Let $x\in R(B)\cap{\frak N}_{\t_1}$. Since $x=y+z,$  $y, z\in
R(B)\cap{\frak N}_{\t_1}$ such that $y^*=y,$ \ $z^*=-z$, then
$$\t_1(ax)=\t_1(ay)+\t_1(az).$$  Thus,  we  need  to  prove  that
$\t_1(E(a/B)x)=\t_1(E(a/B)y)+\t_1(E(a/B)z).$ By Lemma \ref{4.4}, one
finds $\t_1(z)=0$, and due to $(az+za)^*=-(az+za),$ again using
Lemma \ref{4.4} we get $\t_1(az+za)=0,$ which means
$\t_1(az)+\t_1(za)=0$, so $2\t_1(az)=0$, \ i.e. $\t_1(az)=0.$
Consequently, $\t_1(E(a/B)z)=0.$ Since  $y\in B\cap{\frak
N}_{\t_1},$ then by Theorem \ref{4.3} one has $\t_1(ay)=
\t_1(E(a/B)y),$ i.e. $\t_1(ax)=\t_1(E(a/B)x)$ for any $x\in R(B).$
Let $x$ be an arbitrary element taken from   $L_1(W(B),\t_1).$ Then
$x=u+iv$ where $u,v\in L_1(R(B),\t_1)$, so
\begin{eqnarray*}
\t_1(ax)&=&\t_1(a(u+iv))=\t_1(au)+i\t_1(av)\\
&=&\t_1(E(a/B)u)
+i\t_1(E(a/B)v)\\
&=&\t_1(E(a/B)(u+iv))\\
&=&\t_1(E(a/B)x).
\end{eqnarray*}
\end{pf}

Suppose that  $\{A_{\a}\}_{\a\in \mathbb{R}_+}$ is a family of
reversible $JW$- subalgebras of $A$ containing the identity
operator $\id$ such that the set $\bigcup\limits_\a A_{\a}$ is
weak dense in $A$.

\begin{defin}\label{4.6} A family
$\{x_\a\}_{\a\in \mathbb{R}_+}\subset L_1(A,\t)$ is called a
supermartingale if for every $\a \in \mathbb{R}_+$
\begin{enumerate}
\item[(1)] $x_{\a}\in L_1(A_\a,\t)$ \\
\item[(2)] If $\a_1\leq \a_2,$ then $E(x_{\a_2}/A_{\a_1})\leq x_{\a_1}.$
\end{enumerate}\end{defin}

Note that when we replace $\leq$ with the equality sign, then the
family $\{x_\a\}$ is called {\it martingale}. If in Definition
\ref{4.6}, as a $JW$-algebra $A$, we take self-adjoint part of a von
Neumman algebra $W$, then we get usual definitions of
supermartingale and martingale, respectively, in a von Neumann
algebra setting.

It is known \cite{A5} that  $(L_1(W,\t_1))^*=W$, $(L_1(A,
\t))^*=A$, and therefore by $\s_W=\s_W(L_1(W,\t_1), W)$,
$\s_A=\s_A(L_1(A,\t),A)$ we denote weak topologies on
$L_1(W,\t_1)$, $L_1(A,\t)$, respectively. Note that $\s_W$ (resp.
$\s_A$) is generated by a family of seminorms $P_a(x)=|\t_1(ax)|,$
$a\in W$ (resp. $P_a(x)=|\t (a\circ x)|, a\in A$).

\begin{lem}\label{4.7} Let $A$ be a reversible $JW$-algebra with
a f.n.s trace $\t$. Then one has $\s_{W(A)}\mid_{L_1(A,\t)}=\s_A.$
\end{lem}

\begin{pf} Let  $\{x_{\a}\}, \ x\in L_1(A,\t)$ and $x_{\a} \buildrel {\s_{W(A)}} \over
\longrightarrow x$, then $\t_1(ax_\a) \to \t_1(ax)$  for every $a\in
W(A)_+.$ If, in particular, $a\in A_+ \subset W(A)_+$ then
$\t(a\circ x_\a)= {1\over 2} ( \t(a x_\a) + \t(x_\a a) ) = \t(ax_\a)
\to \t(ax)=\t(a\circ x).$ Hence $x_{\a} \buildrel {\s_A} \over
\longrightarrow x$.

Conversely, let $x_{\a} \buildrel {\s_A} \over \longrightarrow x$
and take $b\in W(A)_+$ with $b=c+id$, then  $0\leq c \in A$, \
$d^*=-d$ (see \cite{ARU}). So, $\t(cx_{\a})\to \t(cx)$ since
${dx_{\a}+x_{\a}d \over 2}$ is a skew-symmetric element in $R(A)$
then by Lemma \ref{4.4} one gets
$\t_1(dx_{\a})=\t_1({dx_{\a}+x_{\a}d \over 2})=0$. Hence
$\t_1(bx_{\a}) = \t(cx_{\a})\to \t(cx) = \t_1(bx).$
\end{pf}

Now we are ready to prove a main result of this section.

\begin{thm}\label{4.8}  Let $\{x_\a\}$ be a supermartingale
in $L_1(A, \t).$ If the set $\{x_\a\}$ is weakly relatively compact
in $L_1(A,\t)$, then there is $x\in L_1(A,\t)$ such that
 $x_\a \to x$ in $L_1$ norm.
\end{thm}

\begin{pf} Let $W(A_\a)$ be an enveloping von Neumann algebra
of $A_\a,$ then $x_{\a}\in L_1(A_\a,\t)\subset L_1(W(A_\a),\t_1).$
Moreover, since  $\widetilde E(\cdot/W(A_\a))$ is a conditional
expectation on $W(A)$,  then by Theorem \ref{4.5},  for
$\a_1\leq\a_2$ one has $\widetilde
E(x_{\a_2}/W(A_{\a_1}))=E(x_{\a_2}/ A_{\a_1})\leq x_{\a_1},$ i.e.
$\{x_\a\}$ is a supermartingale in $L_1(W(A),\t_1).$ So Lemma
\ref{4.7} implies that $\{x_\a \}$ is weakly relatively compact in
$L_1(W(A),\t_1).$ Then, Theorem \ref{4.1} yields the existence of
$x\in L_1(W(A),\t_1)$ such that $x_\a \to x$ in $L_1$-norm.
Completeness of $L_1(A,\t)$ w.r.t. $L_1$-norm implies that $x$
belongs to  $L_1(A, \t)$. This completes the proof.
\end{pf}

Note that when the family $\{x_\a\}$ is martingale, the similar
result was proved in \cite{Ber}, therefore, our result extends the
mentioned one. When the trace is finite a similar result was studied
in \cite{A4,A51}

\section*{acknowledgments} The authors thank Prof. V.I. Chilin for his critical reading the paper and given
useful suggestions. The second author (F.M.) acknowledges the MOHE
grant FRGS0308-91. Finally, the authors also would like to thank to
the referee for his useful suggestions which allowed us to improve
the content of the paper.

\end{document}